\documentclass[12pt]{amsart} 

\usepackage[pagebackref=true,colorlinks=true, linkcolor=black ,  citecolor=blue, urlcolor=black]{hyperref}
\usepackage{cleveref} % clever referencing %https://texblog.org/2013/05/06/cleveref-a-clever-way-to-reference-in-latex/
\usepackage[utf8]{inputenc}
\usepackage{bold-extra}

\title{Implementing Real polyhedral homotopy}

\author{Kisun Lee}
\address{School of Mathematical and Statistical Science, Clemson University, 220 Parkway Drive, Clemson, SC 29634, USA}
\email{kisunl@clemson.edu}
\urladdr{https://klee669.github.io}

\author{Julia Lindberg}
\address{Department of Mathematics, University of Texas-Austin, Austin TX USA}
\email {julia.lindberg@math.utexas.edu}\urladdr{https://sites.google.com/view/julialindberg/home}

\author{Jose Israel Rodriguez}\thanks{Research of J.I.R. is supported by the Office of the Vice Chancellor for Research and Graduate Education at U.W. Madison with funding from the Wisconsin Alumni Research Foundation.
}
\address{Department of Mathematics,         University of Wisconsin-Madison,  480 Lincoln Drive, Madison WI 53706-1388, USA.}
\email {jose@math.wisc.edu}\urladdr{https://sites.google.com/wisc.edu/jose}

\usepackage{latexsym,wasysym, mathrsfs, mathtools,hhline,color, bm}
\usepackage[all, cmtip]{xy}

\usepackage{url}

\definecolor{hot}{RGB}{65,105,225}

\usepackage{ textcomp }

\usepackage{graphicx,enumerate}
\renewcommand\emph[1]{\textit{\color{RoyalBlue}#1}}

\usepackage{enumitem}
\usepackage[normalem]{ulem}

 \usepackage{framed} % For leftbar

 \usepackage{amsthm}  
 \usepackage{amsmath} 
 \usepackage{amsfonts}  
 \usepackage{amssymb}
 \usepackage{amscd} 
 \usepackage[all]{xy} 
\usepackage[dvipsnames]{xcolor}  
\usepackage{graphicx}
\usepackage{pdfsync}

\newtheorem{thm}{Theorem}[section]

\theoremstyle{definition}
\newtheorem{ex}[thm]{Example}
\newtheorem{theorem}[thm]{Theorem}

\newtheorem{prop}[thm]{Proposition}

\theoremstyle{definition}
\newtheorem{defn}[thm]{Definition}
\newtheorem{rem}[thm]{Remark}
\numberwithin{equation}{section}

\newcommand{\conv}{{\mathrm{conv}}}
\newcommand{\mvol}{{\mathrm{MVol}}}

\newcommand{\kisun}[1]{{\color{blue}{#1}}}

\newcommand{\R}{\mathbb{R}}

\DeclareMathOperator{\Vol}{Vol}
\DeclareMathOperator{\Log}{Log}
\begin{document}

\maketitle

\begin{abstract}
We implement a real polyhedral homotopy method using three functions.
The first function provides a certificate that our real polyhedral homotopy is applicable to a given system; 
and the second function generates binomial systems for a start system; 
and
the third function outputs target solutions from the start system obtained by the second function. 
This work realizes the theoretical contributions in \cite{ergur2019polyhedral} as easy-to-use functions, allowing for further investigation into real homotopy algorithms.
\end{abstract}

\section{Introduction}\label{sec:introduction}

Finding all points in a zero dimensional algebraic variety is an important problem in many applications in the sciences. This problem amounts to solving a polynomial system of equations with finitely many complex (real or non-real) solutions. 
Many types of algorithms have been proposed using both symbolic and/or numerical techniques. A popular family of numerical algorithms are called \emph{homotopy continuation} algorithms. These algorithms work by continuously deforming the solutions from an `easy' polynomial system into the desired one. While there exists many off-the-shelf homotopy continuation solvers that find all complex solutions, many applications, such as power systems engineering \cite{lindberg2022distribution}, economics \cite{lee2023polyhedral}, and statistics \cite{lindberg2023maximum}, only require knowledge of the real solutions.
In general, there are many more complex solutions than real ones, leading to wasted computation. 
For this reason, the problem of developing 
an efficient homotopy that finds only the real solutions to a polynomial system is an incredibly important open problem in the field of numerical algebraic geometry~\cite{NAG-Book}.

Recent work tackles this problem by presenting an algorithm that certifiably finds all real solutions so long as an %certain
inequality based on the geometry of the polynomial system is satisfied \cite{ergur2019polyhedral}. 
This work relies heavily on  mathematical objects from tropical geometry~\cite{maclagan2021introduction}.
We implement this algorithm in a \texttt{Julia} package, giving the first homotopy based software package that can provably find all real solutions to a patchworked polynomial system without first finding all complex solutions.

We review some of the mathematical concepts behind the algorithms  for the functions in \Cref{sec:background}. In \Cref{sec:rph}, we highlight the key proposition from \cite{ergur2019polyhedral} needed for construction of the \emph{real polyhedral homotopy}. 
In \Cref{sec:functions}, we describe our implementation of the real polyhedral homotopy  which relies on three functions: \texttt{certify\_patchwork}, 
\texttt{generate\_binomials},  
and \texttt{rph\_track}. The code is available as the \texttt{Julia} package \\ \texttt{RealPolyhedralHomotopy.jl}. Documentation for this package can be found at:
\begin{center}
    \url{https://klee669.github.io/RealPolyhedralHomotopy.jl/stable/} 
\end{center}

\section{Preliminaries}\label{sec:background}

\subsection{Regular subdivisions and triangulations}\label{ss:regular-subdivision} 
Let $A\subset\mathbb{Z}^n$ be a set of integer lattice points with convex hull $Q=\conv(A)$. A function $w:A\rightarrow\mathbb{R}$ assigning a real number to each lattice point in $A$ is called a \emph{lifting function}. Denote the collection of lifted points $(a,w(a))\in A\times \mathbb{R}$ by $A^w$. The convex hull of $A^w$ is a polytope, and projecting the lower faces of $\conv(A^w)$ onto the first $n$-coordinates induces a \emph{polyhedral subdivision} $\Delta_w$ of $Q$, i.e.\ all cells in $\Delta_w$ are polyhedral. When all cells of $\Delta_w$ are simplices, it is called a \emph{triangulation}. If a polyhedral subdivision (or triangulation resp.) is induced by a lifting function, we say that the subdivision (or triangulation resp.) is \emph{regular}. 
For a triangulation $\Delta_w$ of $Q=\conv(A)$, we define the \emph{secondary cone} $\mathcal{C}(\Delta_w)$ of $\Delta_w$ as 
the collection of lifting functions that induce the same regular triangulation. Specifically,
\[\mathcal{C}(\Delta_w) = \{ v \in \mathbb{R}^{|A|} \mid \ \Delta_v = \Delta_w \}.\]

\begin{ex}\label{ex:def1}
    
Consider $A = \{0, 1, 2\} \subset\mathbb{Z}$ with lifting function 
    $w(0) = 1, w(1) = 0, w(2) = 3$.  A picture of $A^w$ is given in \Cref{fig:def1}. 
    \begin{figure}[h!]
        \centering
        \includegraphics[width = 0.3\textwidth]{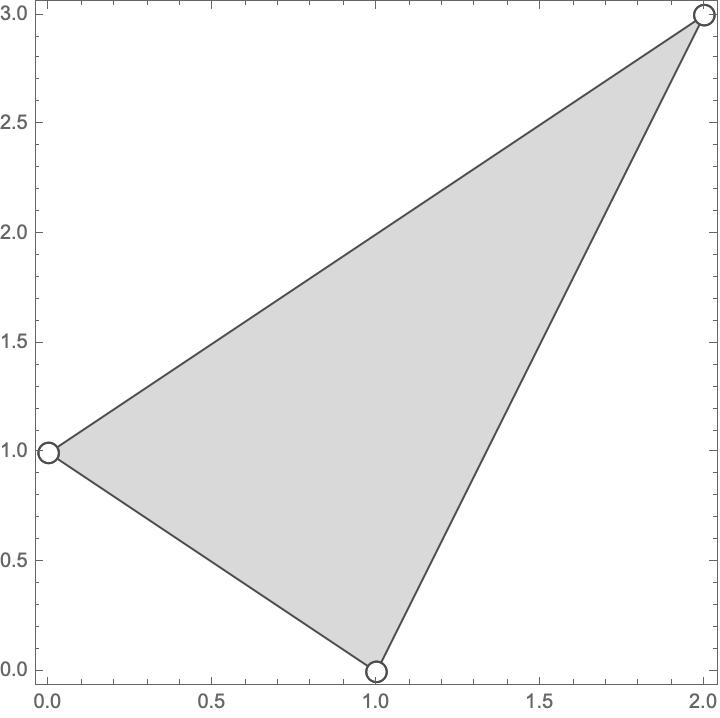}
        \caption{The polytope $A^w$ from \Cref{ex:def1}.}
        \label{fig:def1}
    \end{figure}  
    The lower faces of $\conv(A^w)$ are the line segments $\conv\{(0,1),(1,0)\}$ and $\conv\{(1,0),(2,3)\}$ which induces a regular triangulation of $A$,  $\Delta_w = \{\{0,1\}, \{1,2\}\}$. In this case, we know that any lifting function $w=(w_0,w_1,w_2)$ induces the same triangulation if $w_0\geq w_1$ and $w_2\geq w_1$ when $w_0,w_1,w_2$ are not all equal.
    Therefore, we have
    \[\mathcal{C}(\Delta_w) = \{(w_0,w_1,w_2) \in \mathbb{R}^3 \ | \ w_0 \geq w_1, w_2 \geq w_1, \ w_0, w_1, w_2 \ \text{ are not all equal} \}. \]
\end{ex}

\subsection{Cayley configurations and mixed cells}\label{ss:Cayley-configuration} 
For sets of lattice points $A_1,\dots, A_m$ in $\mathbb{Z}^n$, consider \kisun{the} set 
\[A_1*\cdots*A_m := \{(a_i,e_{n+i})\in \mathbb{Z}^{n+m}\mid a_i\in A_i, i=1,\dots,m\}\]
of lattice points in $\mathbb{Z}^{n+m}$
where $e_{n+i}$ is the $(n+i)$-th canonical vector. The set $A=A_1*\cdots *A_m$ is called a \emph{Cayley configuration}. 
Given a Cayley configuration $A \in \mathbb{Z}^{n+m}$, we lift $A$ with a lifting function $w$ and construct a polyhedral subdivision $\Delta_w$ of $\text{conv}(A)$. A cell $\sigma$ in $\Delta_w$ is called \emph{mixed} if $\sigma$ has exactly two elements from each $A_i$. 
For a fixed cell $\sigma$, we are interested in all lifting functions $w$ such that $\sigma$ is a mixed cell of $\Delta_w$.
We call such set a \emph{mixed cell cone of $\sigma$}, and it is formally defined by
\[M(\sigma):= \{w\in \mathbb{R}^{|A|}\mid \sigma \text{ is a mixed cell in }\Delta_w\}.\]
For a triangulation $\Delta_w$ of $A$ induced by $w$, we define the \emph{mixed cell cone of $\Delta_w$} as 
\[M(\Delta_w):=\bigcap\limits_{\{\sigma \text{ is a mixed cell in }\Delta_w\}}M(\sigma).\]

\begin{ex}\label{ex:def2}
    Consider $A = \{0,1,2\}$ with lifting function $w = (1,0,3)$ as in \Cref{ex:def1}. Since $m = n = 1$, 
    the Cayley configuration of $A$ is just a translation of $A$. We observe that $\Delta_w$ is mixed since each cell contains exactly two elements from $A$. 
    Let $\sigma_1 = \{0,1\} \in \Delta_w$ and $\sigma_2 = \{1,2\} \in \Delta_w$. Then $M(\sigma_1) = M(\sigma_2)$, giving
    \[M(\Delta_w) = \{(w_0,w_1,w_2) \in \mathbb{R}^3 \ | \ w_0 \geq w_1, w_2 \geq w_1, \ w_0, w_1, w_2 \ \text{ are not all equal} \} \]
    which is the same as $\mathcal{C}(\Delta_w)$.
\end{ex}

\subsection{\texorpdfstring{$A$}--discriminant and its amoeba}\label{ss:A-disc} 
For a sparse polynomial $f\in \mathbb{C}[x_1,\dots, x_n]$, we define the (monomial) \emph{support} $A_f$ of $f$ as %which is
the set of exponents of all monomials of $f$. Let $\mathbb{C}^{|A_f|}$ be the set of polynomials with complex coefficients supported on $A_f$.
Then, a polynomial $f\in \mathbb{C}^{|A_f|}$ can be written as
\[f(x) = \sum\limits_{a\in A_f}c_ax^a,\]
where $x^a$ is the monomial $x_1^{a_1}\cdots x_n^{a_n}$ for $a=(a_1,\dots,a_n)$.
For a polynomial $f=\sum\limits_{a\in A_f}c_ax^a$ supported on $A_{f}$, let $c_f:=(c_a)_{a\in A_f}$
be the coefficient vector for $f$. 
For an $n$-tuple of vectors $C=(c_{f_1},\dots, c_{f_n})$, we let $F_C:=\langle f_1,\dots, f_n\rangle$ be a square polynomial system such that $c_{f_i}$ is a coefficient vector for $f_i$ for each $i=1,\dots,n$. 
For a Cayley configuration $A = A_{f_1}*\cdots* A_{f_n}$, we  define the \emph{$A$-discriminant} 
\[\nabla_{A}:=\left\{C\in \mathbb{C}^{|A|} \mid \ F_C(x) \text{ has a singularity in }(\mathbb{C}\setminus\{0\})^n\right\}.\]

We say that $\nabla_A$ is \emph{non-defective} if it is codimension one.
We are particularly interested in $\nabla_A(\R):=\nabla_A\cap \R^{|A|}$ 
because in a connected component of the complement of $\nabla_A(\R)$, 
the number of real solutions to the corresponding polynomial systems is constant.

When an $A$-discriminant is non-defective, we may consider the Newton polytope of its
defining polynomial. For a vertex $v$ of the Newton polytope, its normal cone is denoted by $\mathcal{N}(v)$. For a lifting function $w$ for $A$, 
we note the following relations, proved in \cite[Lemma 2.16]{ergur2019polyhedral}, between the secondary cone $\mathcal{C}(\Delta_w)$, mixed cell cone $M(\Delta_w)$ and the normal cone $\mathcal{N}(v)$ of $\nabla_A$:
\begin{equation}\label{eq:containments}
\mathcal{C}(\Delta_w)\subseteq M(\Delta_w)\subseteq\mathcal{N}(v)
\end{equation}
for the vertex $v$ of the Newton polytope of $\nabla_A$ satisfying that $w\in\mathcal{N}(v)$. 

The \emph{log-absolute value map} $\Log|\cdot| :(\mathbb{C}\setminus\{0\})^n\rightarrow \mathbb{R}^n$ is defined by
\begin{align}
\Log | (x_1,\dots, x_n)|:=(\log|x_1|,\dots, \log|x_n|).\label{def:log_lift}
\end{align}
For a Laurent polynomial $f\in \mathbb{C}[x_1^\pm,\dots, x_n^\pm]$, the image of the variety $\mathcal{V}(f)$ under the log-absolute value map is called the \emph{amoeba} of $f$ and denoted by $\mathcal{A}(f)$. For a textbook reference, see \cite[Chapter 1]{maclagan2021introduction}.%. 

We remark that the complement of an amoeba consists of convex regions. For the real polyhedral homotopy algorithm, we will consider a specific lifting function $w$ inducing a triangulation $\Delta_w$ such that mixed cell cone $M(\Delta_w)$ is contained in a complement of the amoeba of $\nabla_A$.

\subsection{Homotopy continuation} For a system $F=\langle f_1,\dots, f_n\rangle$ of polynomials in $n$ variables, finding all isolated solutions of the system is an important task.
\emph{Homotopy continuation} is a method 
to find numerical approximations of solutions of a system of polynomial equations. The main idea is to track solution paths from a system $G$ called a \emph{start system} whose solutions are known to the \emph{target system} $F$. For a start system $G=\langle g_1,\dots, g_n\rangle$ with the same number of variables and equations of $F$, we construct a homotopy $H(x,t)$ such that $H(x,0)=G$ and $H(x,1)=F$. 
To track the solutions from $G = 0$ as $t$ varies from zero to one, a common approach is to use \emph{predictor-corrector} methods. These methods rely on numerically solving an ordinary differential equation, called the Davidenko equation, as well as using Newton iterations.
For details, see \cite[Chapter 2]{SommeseWampler:2005}. 

\subsection{Polyhedral homotopy continuation and Bernstein's theorem}\label{ss:phc} 
In order to choose a start system $G$ for a homotopy continuation algorithm, we want the number of solutions of $G(x) = 0$ to be roughly equal to the number of solutions of $F(x) = 0$.
In this paper, the \emph{polyhedral homotopy continuation} established by Huber and Sturmfels \cite{huber1995polyhedral} is considered. 

For polytopes $P,Q$ in $\mathbb{R}^n$, the \emph{Minkowski sum} of the polytopes is defined as
$P+Q=\{p + q \mid p\in P, \ q \in Q\}$. 
For a nonnegative number $a$ and polytope $P$, we define
$aP:= \{a p \mid p\in P \}$.
The Euclidean volume of the Minkowski sum $\Vol(a_1Q_1+\cdots +a_nQ_n)$ is a homogeneous polynomial  in variables $a_1,\dots, a_n$. The coefficient of the mixed term $a_1a_2\cdots a_n$ of the polynomial $\Vol(a_1Q_1+\cdots +a_nQ_n)$ is called the \emph{mixed volume} of $Q_1,\dots, Q_n$ and is denoted by $\text{MVol}(Q_1,\dots, Q_n)$. The following celebrated theorem 
relates the mixed volume of the Newton polytopes of a polynomial system to the number of isolated solutions in the torus.
\begin{theorem}[Bernstein's theorem]\cite[Theorem A]{bernshtein1975number}
  Let $F$ be a system of polynomials $f_1,\dots, f_n$ in $\mathbb{C}[x_1,\dots, x_n]$.
  For Newton polytopes $Q_{f_i}$ for each $f_i$, we have 
  \[(\text{the number of isolated solutions of }F\text{ in }(\mathbb{C}\setminus\{0\})^n)\leq \text{MVol}(Q_{f_1},\dots, Q_{f_n}).\]
  Furthermore, for polynomials $f_1,\dots, f_n$ with generic coefficients the inequality
is tight.
\end{theorem}

Polyhedral homotopy continuation tracks $\mvol(Q_{f_1},\dots, Q_{f_n})$  paths to find all solutions to $F = \langle f_1,\ldots, f_n \rangle = 0$.
The idea is to construct 
a collection of binomial start systems whose number of solutions sum to  $\mvol(Q_{f_1},\dots, Q_{f_n})$. We now describe how to find one such system, $G$.
Consider a polynomial 
\[f(x) = \sum\limits_{a\in A}c_ax^a\]
supported on $A$.
Consider a lifting function $w$. By multiplying each monomial of $f$ by $t^{w(a)}$,
we have the  \emph{lifted polynomial}
\[\overline{f}(x,t)=\sum\limits_{a\in A}c_ax^a t^{w(a)}\]
Suppose that a square polynomial system $F$ consists of polynomials $f_1,\dots, f_n$ supported on $A_{f_1},\dots, A_{f_n}$, respectively. 
Lifting all polynomials 
in $F$ gives the lifted system $\overline{F}(x,t)$ which satisfies $\overline{F}(x,1)=F$. 
The solutions of $\overline{F}$ can be expressed by $x(t)=(x_1(t),\dots, x_n(t))$ where each $x_i(t)$ is a Puiseux series and
\[x_i(t)=t^{\alpha_i}y_i+ \text{(higher order terms in $t$)}\]
for some rational number $\alpha_i$ and a nonzero constant $y_i$. Plugging $x(t)$ into polynomials gives
\[\overline{f}_j(x(t),t) =\sum\limits_{a\in A_{f_j}}c_a y^a t^{\langle a,\alpha\rangle+w(a)}+\text{(higher order terms in $t$}).\] where $y^a=y_1^{a_1}\cdots y_n^{a_n}$.
For the minimum value of $\langle a,\alpha\rangle +w(a)$ over all $a\in A_{f_j}$, divide by $t^{\langle a,\alpha\rangle +w(a)}$ and set $t=0$. Iterating this for each $\overline{f}_j$, we have a start system~$G$. Note that for a fixed lifting function $w$, the minimum can be obtained from different monomials $a\in A_{f_j}$ depending on $\alpha$. Therefore, we may have several start systems $G$ and they induce a collection of start systems.

For most choices of $\alpha$, the procedure outlined above gives a start system consisting of monomials, which is not useful since monomial systems of equations have no solutions in the torus. Instead, one chooses $\alpha$ carefully to get a binomial start system since general binomial systems of polynomial equations have solutions in the torus and can be solved efficiently using linear algebra~\cite{chen2014solutions}.
Such an $\alpha$ and the corresponding binomial system can be obtained from a mixed cell of a triangulation $\Delta_w$ of $\conv(A)$ induced by $w$.

We summarize this section with the following definition. 

\begin{defn}[Polyhedral homotopy]\label{def:ph}
For a polynomial system $F=\langle f_1,\dots,f_n\rangle$, with $f_i=  \sum\limits_{a\in A_{f_i}} c_{a}x^{a}$,
define the \emph{polyhedral homotopy} as 
\begin{align*}
    H_m(x,t) := \begin{cases}
    \sum\limits_{a\in A_{f_1}} t^{m_{1,a}}c_{a}x^{a} \\
    \qquad \vdots \vspace{.3pc}\\
     \sum\limits_{a\in A_{f_n}} t^{m_{n,a}}c_{a}x^{a}.
    \end{cases}
\end{align*}
where each element of $m= \{  m_{i,a} \mid i\in [n]\kisun{,}\text{ and } a\in A_{f_i}\}$ is in $\mathbb{Q}_{\geq 0}$ and for each $i\in [n]$, the minimum of $\{ m_{i,a} \mid a\in A_{f_i}\}$
equals zero and is obtained precisely twice.
We call $m$ the \emph{lifting of the polyhedral homotopy}. 
\end{defn}

By substitution, we have $H_m(x, 1) \equiv F (x)$, which is the target system. Assuming genericity of the lifting $m$, and system $F$, $H_m(x,t)$ has smooth, non-intersecting paths in $(\mathbb{C}\setminus\{0\})^n \times  (0, 1]$,  parameterized by $t$.
The starting points of these paths at $t=0$ can be obtained by solving binomial systems. The polyhedral homotopy continuation is implemented in software  \texttt{HOM4PS2} \cite{lee2008hom4ps}, \texttt{PHCPack} \cite{verschelde1999algorithm}, \texttt{Pss5}~\cite{Pss5}, and \texttt{HomotopyContinuation.jl}~\cite{breiding2018homotopycontinuation}.

\subsection{Certification of numerical solutions}\label{sec:certification}

Practical implementations of homotopy continuation 
rely on heuristics to compute  
numerical approximations of solutions to a system of polynomial equations.
Therefore, \textit{a posteriori} certification for the correctness of these approximations may be necessary. 
We say a numerical approximation of a solution of a system $F$ is \emph{certified} if it can be refined up to arbitrary precision to an actual solution of $F$ that is uniquely contained in a neighborhood of the approximation by applying a suitable operator (e.g., the Newton operator). Smale's $\alpha$-theory \cite[Chapter 8]{blum2012complexity} and the Krawczyk method \cite[Chapter 8]{moore2009introduction} are commonly used for certification. Numerical root certification algorithms are implemented in software \texttt{alphaCertified} \cite{hauenstein2011alphacertified}, \texttt{NumericalCertification.m2} \cite{lee2019certifying} and the function \texttt{certify} in \texttt{HomotopyContinuation.jl} \cite{breiding2020certifying}.

\section{Finding a real polyhedral homotopy}\label{sec:rph}

The real polyhedral homotopy algorithm introduced in \cite{ergur2019polyhedral} provides a framework for finding all real solutions of a polynomial system devised as a variation of the standard polyhedral homotopy. The 
main idea stems from Viro's patchworking for complete intersections \cite{sturmfels1994viro}. 
This result establishes  
a homeomorphic relation between the set of real solutions of a polynomial system and the intersection of the positive part of tropical varieties (see \cite[Section 2.2]{ergur2019polyhedral} for a rigorous statement).
We say a system is \emph{patchworked} if the real solution set  of the system is homeomorphic to the union of the intersection points of the subcomplexes
obtained by Viro's patchworking.
In other words, we are interested in constructing a homotopy whose number of real solutions does not change while tracking.

Sections \ref{ss:regular-subdivision}, \ref{ss:Cayley-configuration} and \ref{ss:A-disc},
provide the preliminaries for the following proposition motivating the real polyhedral homotopy algorithm in \cite{ergur2019polyhedral}.

\begin{prop}\cite[Proposition 2.19]{ergur2019polyhedral}\label{prop:frameworkRPH}
Consider a polynomial system $F_C=\langle f_1,\dots, f_n\rangle$ with a coefficient vector $C$ and support sets $A_{f_1},\dots, A_{f_n}$ such that $\dim(A_{f_i})=n$ for each $i$. Suppose that $u=(u_1,\dots, u_n)$ where $u_i=(u_a)_{a\in A_{f_i}}$ is a vector satisfying that
\begin{enumerate}
    \item The vector $u$ is not on the boundary of any secondary cone of the Cayley configuration $A=A_{f_1}*\cdots*A_{f_n}$.
    \item The ray $\Log|C|+\lambda u$ does not intersect $\mathcal{A}(\nabla_A(\R))$ for any $\lambda\in[0,\infty)$.
\end{enumerate}
Then, the tuple of real Puiseux series $x(t)=(x_1(t),\dots, x_n(t))$ where $x_i(t)$ is a real Puiseux series and
\[x_i(t)=t^{\alpha_i}y_i +(\text{higher order terms in $t$})\]
is a solution to the system $F_C$ only if $(\alpha_1,\dots, \alpha_n, 1)$ is an outer normal to a lower facet of $\sum_{i=1}^n \text{conv}(A_i^{u_i})$.
\end{prop}
The solution $x(t)$ in \Cref{prop:frameworkRPH} is derived in a similar way to the polyhedral homotopy method outlined in \Cref{ss:phc}.

Since the conditions from the proposition above are not satisfied in general, a certification procedure for checking if a polynomial system is patchworked is suggested in \cite{ergur2019polyhedral}.
Note that the containment (\ref{eq:containments}) shows that the mixed cell cone $M(\Delta_w)$ is contained in the corresponding normal cone $\mathcal{N}(v)$ of $\nabla_A$. Therefore,
a proper choice of a lifting function $w$ might result in the mixed cell cone $M(\Delta_w)$ that is contained in the complement of $\mathcal{A}(\nabla_A(\mathbb{R}))$.
When such a coefficient vector is given,
costly computation of the amoeba $\mathcal{A}(\nabla_A(\mathbb{R}))$ can be replaced by a mixed cell cone computation.
Based on the argument above, the following proposition establishes a sufficient condition for a polynomial system to be patchworked.

\begin{prop}\cite[Proposition 3.1]{ergur2019polyhedral}\label{prop:certification}
Let $F_C=\langle f_1,\dots, f_n\rangle$ be a system of sparse polynomials with coefficient vector $C$ with support sets $A_{f_1},\dots, A_{f_n}$. Let $\Delta_w$ be the triangulation of the Cayley configuration $A=A_{f_1}*\cdots *A_{f_n}$ induced by the lifting $w=\text{Log}|C|$. Consider the corresponding dual mixed cell cone $M(\Delta_w)^\circ$ and its generating vectors $\zeta_1,\dots, \zeta_L$. Then, 
\begin{align}\label{eq:log_ineq}
    \langle\Log|C|,\zeta_i\rangle>\|\zeta_i\|_1\log(|A|) 
\end{align}
for all $i=1,\dots, L$ implies that the system $F_C$ is a patchworked system. Also, for any $v\in M(\Delta_w)$, the ray $\text{Log}|C|+\lambda v$ for $\lambda\in [0,\infty)$ does not intersect $\mathcal{A}(\nabla_A(\mathbb{R}))$.
\end{prop}

\begin{defn}[Real polyhedral homotopy]
We say a lifting $m$ induces a \emph{real polyhedral homotopy} of $F = \langle f_1,\ldots, f_n \rangle$ if
\begin{itemize}
    \item   $H_m(x, 1) \equiv F (x)\subset \mathbb{R}[x]$
    \item $ \{ (x,t) \in (\mathbb{R}\setminus\{0\})^n \times (0,1] \mid  H_m (x, t) = 0, \text{ and } t\in (0,1]\}$ defines smooth non-intersecting paths in $(\mathbb{R}\setminus\{0\})^m \times  (0, 1]$, 
    parameterized by $t$, and emanating from isolated $(\mathbb{R}\setminus\{0\})$-zeros of $F(x) \equiv H(x,1)$ and continue toward $t = 0$. 
    \item the starting points of these paths as $t=0$ are obtained by solving  binomial systems coming from mixed cells and the $(\mathbb{R}\setminus\{0\})$-zeros of the target system $F(x)$ can be found by tracking these paths over the real numbers.
\end{itemize}
\end{defn}
Proposition~\ref{prop:certification} gives a method to certify when the lifting $m = \log(|C|)$ induces a real polyhedral homotopy.

\begin{ex}
If $f_1 = -1 - 24000y + x^3$,
$f_2 = -9 + 50xy - y^2$, 
then the $\Log |C|$ lifting corresponds to the following lifted polynomials:
\begin{align*}
- 1 - 24000t^{\log(24000)}y + x^3,  \\
 -9 t^{\log(9)}  + 50t^{\log(50)} xy -y^2.
\end{align*}
Denoting the real Puiseux series solutions by $(x(t), y(t))$, these lifted polynomials give two homotopies with binomial start systems induced by mixed cells. The homotopies are:
\begin{align*}
h_1(x,y ; t) &= \langle- t^{a_1} - 24000y +x^3,  
 -9 t^{a_2}  + 50 xy -y^2\rangle \\
h_2(x,y ;t) &= \langle- t^{b_1} - 24000y +x^3, -9   + 50 xy -t^{b_2} y^2\rangle
\end{align*}
where $a_1,a_2,b_1,b_2 \in \mathbb{Z}_{>0}$. We implement a way to certify that this is a (well-chosen) real polyhedral homotopy and outline our implementation in the next section. 
\end{ex}

\section{Real polyhedral homotopy implementation}\label{sec:functions}

This section describes our \texttt{Julia} implementation in detail. We assume we are given a polynomial system $F = \langle f_1,\ldots, f_n \rangle $ with finitely many complex solutions. 
Let $C$ denote the vector of coefficients of the system $F$.

\subsection{\texttt{certify\_patchwork} }
This function certifies 
if a given system is patchworked so that all real solutions can be found using the real polyhedral homotopy. It checks if inequality \eqref{eq:log_ineq} holds for each mixed cell.
It returns the value $1$ if the system $F$ is certified to be patchworked according to the inequality. 
Otherwise, $0$ is returned.

\begin{leftbar}
	\begin{verbatim}
@var x y
f1 = -1 - 24000*y + x^3
f2 = -9 + 50*x*y - y^2
F = System([f1, f2])
certify_patchwork(F)
    1
	\end{verbatim}
\end{leftbar}

\begin{rem}
\label{rmk:Number_Real_Solutions}
As an optional argument we have \texttt{Number\_Real\_Solutions}. 
When this optional argument is set to true (default is false) 
we return $(1,k)$ where $k$ is number of real solutions to the target system when the target system is patchedworked.
 It works by solving the binomial systems $B$ (discussed in  \Cref{ss:generate_binomials})  by  using Smith normal forms as outlined in \cite[Section 2.5]{ergur2019polyhedral}.   For additional details on solving binomial systems see \cite{chen2014solutions}.
 Otherwise, we return $0$. 
 \begin{leftbar}
	\begin{verbatim}
certify_patchwork(F; Number_Real_Solutions = true)
    (1,4)
	\end{verbatim}
\end{leftbar}

\end{rem}

\subsection{\texttt{generate\_binomials}}\label{ss:generate_binomials}
This function takes as an input a polynomial system $F = \langle f_1,\ldots, f_n \rangle$ in $n$ variables and outputs an object called \texttt{Binomial\_system\_data} that consists of 4 objects: a collection of binomial systems, their normal vectors, the lifting function, and the mixed cells. 
Contents in the object \texttt{Binomial\_system\_data} stem from the mixed cells induced by the $\Log|C|$-lifting mentioned in \Cref{prop:certification}. Each value in \texttt{Binomial\_system\_data} can be called by its name \texttt{binomial\_system}, \texttt{normal\_vectors}, \texttt{lifts}, and \texttt{cells} as follows:

\begin{leftbar}
	\begin{verbatim}
# Continued from above.
B = generate_binomials(F);
B.binomial_system
    2-element Vector{Any}:
    Expression[-24000*y + x^3, 50*x*y - y^2]
    Expression[-24000*y + x^3, -9 + 50*x*y]	    
B.lifts
    2-element Vector{Vector{Int64}}:
    [0, -10085809, 0]
    [-3912023, 0, -2197225]    \end{verbatim}
\end{leftbar}

In our implementation, we use the \texttt{Julia} package \texttt{MixedSubdivisions.jl} \cite{mixed-subdivisions} to compute the mixed cells. 
The package requires a lifting function to take on integer values.
So in our implementation, we take 
$10^6 \cdot \Log|C| $ rounded to the nearest integer as our lifting. 
Using this scaled and rounded lifting is suitable because uniformly scaling a lifting function gives the same triangulation.
Therefore, rounding the lifting function will preserve the desired triangulation so long as $\Log|C|$ is of distance at least $10^{-6}$ from the boundary of its mixed cell cone. We emphasize that this is a heuristic and for some systems, one may need to take $10^k \cdot \Log|C|$ for $k>6$. In all of our experiments, $k = 6$ was sufficient to preserve the desired triangulation.

\subsection{\texttt{rph\_track}}
Our function \texttt{rph\_track} takes an object \texttt{Binomial\_system\_data} $B$ and a polynomial system $F$ as its input. It returns the output of tracking the real solutions of the binomial systems to the target system. For our running example, the function returns all real solutions to $F=0$ and no nonreal solutions.

\begin{leftbar}
	\begin{verbatim}
rph_track(B,F)
    4-element Vector{Vector{Float64}}:
    [-1095.4451129504978, -54772.25548320812]
    [1095.4451137838312, 54772.255524874796]
    [8.111114476617955, 0.02219298606763958]
    [-8.103507635567631, -0.022213821121964985]    
	\end{verbatim}
\end{leftbar}

In contrast, the default solve command in \texttt{HomotopyContinuation.jl} finds all six complex solutions of which two are nonreal and the four real solutions are the solutions found above. 
\begin{leftbar}
	\begin{verbatim}
S = HomotopyContinuation.solve(F)
solutions(S)
    6-element Vector{Vector{ComplexF64}}:
    [-0.003803837191831332 - 8.10709081578636im,
                        -1.0415806358129898e-5 + 0.022201564813120536im]
    [-0.003803837191831332 + 8.10709081578636im,
                        -1.0415806358129898e-5 - 0.022201564813120536im]
    ...
	\end{verbatim}
\end{leftbar}

We implement our function in two steps:
The first step is to use Smith normal forms to find all real solutions to the binomial systems as mentioned in \Cref{rmk:Number_Real_Solutions}. 
The second step of this function tracks the real solutions from the binomial systems to the target system using \texttt{HomotopyContinuation.jl}.

\begin{rem}
In the event \texttt{certify\_patchwork} function returns $0$, \texttt{rph\_track} function will still run. In this situation, there is no guarantee that the real solutions of the start system will converge to real solutions of the target system nor that every real solution of the target system will be returned. 
The following example shows how the real polyhedral homotopy behaves for a non-patchworked system.
\begin{leftbar}
	\begin{verbatim}
@var x y
F2 = System([ -1 - 240*y + x^3,  -9 + 50*x*y - y^2 ]
certify_patchwork(F2)
    0
    
B2 = generate_binomials(F2);
R = rph_track(B2,F2)
    4-element Vector{Vector{Float64}}:
    [-109.54445340262997, -5477.221026962461]
    [109.54453673601331, 5477.225193632879]
    [2.601807483849416, 0.06921950525397731]
    [-2.525776652013741, -0.07130546925252307]
\end{verbatim}
\end{leftbar}
As a consequence, for non-patchworked system we can view our implementation as a heuristic for finding real solutions to a polynomial system. 
\end{rem}

\subsection{An \texttt{rph\_track} option}
The optional argument \texttt{Certification}, certifies all real solutions to a patchworked system $F$ found through the aforementioned real polyhedral homotopy algorithm are in fact good approximations to $F(x) = 0$.  
Given correct path tracking, the real polyhedral homotopy guarantees to find all real solutions when the system is patchworked.
To guarantee that the path tracking was done correctly, we use an \textit{a posteriori} certification.
For a patchworked system $F$, when the real polyhedral homotopy root-finding is certified, the function returns a list of solutions to $F(x) = 0$ and $1$; otherwise, it returns $0$. 

\begin{leftbar}
	\begin{verbatim}
rph_track(B,F; Certification = true)
    ([[-1095.4451129504978, -54772.25548320812],...], 1)
	\end{verbatim}
\end{leftbar}

When the option \texttt{Certification} is set to \texttt{true}, two extra steps are done in $\texttt{rph\_track}$. 
First, it certifies the real start solutions of a binomial system in $B$ are true solutions to $B$, which provides the real root count of $F=0$.
For the second additional step, the numerical approximations for solutions \kisun{of} $F$ produced by the real polyhedral homotopy tracking are certified using the Krawczyk method.
The certification relies on  \texttt{certify} in \texttt{HomotopyContinuation.jl} with details 
found in \cite{breiding2020certifying}.

\subsection{Reliability}\label{ss:reliability}
For those who are not frequent users of numerical methods in algebraic geometry, we add this subsection to clearly explain the reliability of our package and the meaning of certification of the output. 
Our package relies on the path-tracking heuristics in \texttt{HomotopyContinuation.jl}. 
It is well known that numerical path-tracking methods sometimes miss solutions. While in theory such an event happens with probability zero,  in practice, we work with floating point arithmetic so may observe this behavior.  

Both verifying (like a trace test \cite{leykin2018trace,brysiewicz2023sparse}) and certifying the completeness of the set of real solutions still remain open, to the best of the authors' knowledge. 
The certification functionality in the package is limited to checking if all numerical solutions that are computed are certified in the context of \Cref{sec:certification}. 
When the system is not patchworked, our implementation does not guarantee finding all real solutions to a given system, even if path-tracking was perfect. Nevertheless, our package can still run the real polyhedral homotopy and get an output that is a (potentially proper) subset  of the real solution set. When the certification option is turned on the output consists of certified real solutions with no guarantee of completeness of the real solution set.

\section{Outlook}

There may be potential improvement available for the real polyhedral homotopy algorithm described in this paper by constructing homotopy paths that avoid the real discriminant locus. The current  verification using inequality \eqref{eq:log_ineq} 
does not give both necessary and sufficient conditions for detecting if a polynomial system is patchworked and as a result is a conservative way to certify this.
To illuminate the potential of the real polyhedral homotopy algorithm described in this paper, we provide an example of a polynomial system from game theory where the real polyhedral homotopy algorithm finds all real solutions even though the system is not certified to be patchworked.

\begin{ex}[A modified version of \cite{nash2023}]\label{ex:nash}
    Consider a three-player game such that each player has $3$ strategies. For $i=1,2,3$, define the $i$-th player's payoff by a $3\times 3\times 3$-tensor $P^{(i)}$, and let $p^{(i)}_j$ be the probability that the $i$-th player chooses the $j$-th strategy. In this case, the $i$-th player's payoff $a_i$ can be computed as:
    \[a_i=\sum_{j,k,l=1}^3P^{(i)}_{j,k,l}\cdot p^{(1)}_jp^{(2)}_kp^{(3)}_l\quad\text{for}\quad i=1,2,3.\]
    Assuming that all players do not change their initial strategy, the vector of probabilities for strategies of players optimizing output of their own is called a \emph{Nash equilibrium}. The problem of finding all Nash equilibria can be modelled as a constrained optimization problem. Constrained optimization problems with smooth critical points can be solved using the method of Lagrange multipliers. This method then reduces to solving a system of polynomial equations. For example, all Nash equilibria for the three-player game described above appear as real solutions to the following system of $12$ polyomial equations in $12$ unknowns $(p^{(1)}_1,p^{(1)}_2,p^{(1)}_3,\,\,
    p^{(2)}_1,p^{(2)}_2,p^{(2)}_3,\,\,
    p^{(3)}_1,p^{(3)}_2,p^{(3)}_3,\,\,
    a_1,a_2,a_3)$:
    \begin{align*}
        \sum_{j=1}^3 p^{(1)}_j=\sum_{j=1}^3p^{(2)}_j=\sum_{j=1}^3p^{(3)}_j=1&\\
        p^{(1)}_j\left(a_1-\sum_{k,l=1}^3P^{(1)}_{j,k,l}\cdot p_k^{(2)}p_l^{(3)}\right)=0&\quad\text{for }j=1,2,3\\
        p^{(2)}_k\left(a_2-\sum_{j,l=1}^3P^{(2)}_{j,k,l}\cdot p_j^{(1)}p_l^{(3)}\right)=0&\quad\text{for }k=1,2,3\\
        p^{(3)}_l\left(a_3-\sum_{j,k=1}^3P^{(3)}_{j,k,l}\cdot p_j^{(1)}p_k^{(2)}\right)=0&\quad\text{for }l=1,2,3.         
    \end{align*}
For randomly chosen payoff tensors
\begin{align*}
    P^{(1)}&=\left[\begin{array}{ccc|ccc|ccc}
         34 &   162 &  70 &     140   &  174 & -183 &   104   & 68 & 18  \\
     136 & 2001 & 72 &  166 & 140  &  -10541 &  1261  & 72  &20 \\
      80  &  32  &  10 &     132   &  176  & 17 &  174 & 86 & 94
     \end{array}\right],\\
         P^{(2)}&=\left[\begin{array}{ccc|ccc|ccc}
  56 &-6& 22& 5712& -150& 74 &150& 10& 150\\
  -28& 80& 1196& 182& 116& 22& 44& 186& 60\\
  12& 48& 22& 64& 11& 46 &80& 1192& 82
     \end{array}\right],\\
         P^{(3)}&=\left[\begin{array}{ccc|ccc|ccc}
  104 &162& 655& 2116& 130& 559& 134& 138& -4\\
  172& 152& -12& 25& 188& 162& 1168& 162& 178\\
  134& 124& 170& 118& 38& 130& 113& 800& 152
     \end{array}\right],
\end{align*}
we solve this system using \texttt{RealPolyhedralHomotopy.jl} and see that it returns all eight real solutions.
\begin{leftbar}
    \begin{verbatim}
8-element Vector{Vector{Float64}}:
 [0.023564140817086864, -0.003578459398610543,
                        ... , 17.00815164148862, 150.9733200702087]
 [0.1690458198092814, 0.94412282597982,
                        ... , 40.84440318232266, 199.80455930749537]
 ...
    \end{verbatim}
\end{leftbar}
In this case, the first nine coordinates of solutions to the polynomial system correspond to probabilities. By restricting the real solutions to those with positive first nine coordinates we find one Nash equilibrium.
\begin{leftbar}
    \begin{verbatim}
valid_real_solutions = filter(s -> all(s[1:9] .> 0), rsols)
    [0.38164510348087044, 0.3858370607366621, 
                        ... , 221.7388487672952, 210.9466957898507]
\end{verbatim}
\end{leftbar}
\end{ex}
The success of Example~\ref{ex:nash} motivates the future development of heuristic methods for finding a complete real solution set. In addition to applications in economics, as seen in the previous example, we hope that problems in reaction networks, statistics, and optimization, where polyhedral geometry has played a role in counting complex solutions, are ripe for applying real polyhedral homotopy.

\section*{Acknowledgment}

The authors are grateful to Paul Breiding, Timo de Wolff and Christopher O'Neill for helpful discussions.

\bibliographystyle{abbrv}
\bibliography{refs.bib}
\end{document}